\documentclass[twoside,bezier]{amsart}
\usepackage[left=2cm,right=2cm]{geometry}

\usepackage{amsmath,amsthm,amsfonts,amssymb,amscd,euscript}
\usepackage{setspace}

\usepackage[figuresright]{rotating}
\newtheorem{thm}{Theorem}[section]

\newtheorem{lem}[thm]{Lemma}

\newtheorem{exam}[thm]{Example}
\theoremstyle{definition}\newtheorem{defn}[thm]{Definition}
\theoremstyle{remark}
\newtheorem{rem}[thm]{Remark}
\numberwithin{equation}{section}

\newcommand{\A}{\mathcal{A}}

\newcommand{\filebegin}{\begin{document}}
\newcommand{\fileend}{\end{document}}
\makeatother
\def\thefootnote{}
\newcommand{\lo}{\longrightarrow}
\newcommand{\NMM}{\hspace*{2mm}}

\usepackage[left=2cm,right=2cm]{geometry}

\def\n{\noindent}%

\numberwithin{equation}{section}
\def\mapdown#1{\Big\downarrow\rlap
{$\vcenter{\hbox{$\scriptstyle#1$}}$}}
\begin{document}

\vspace*{2cm}
\begin{center}
{\bf\large
Modular cone metric spaces}
 \\[0.5cm]
{Saeedeh Shamsi Gamchi and Mohammad Janfada, Asadollah Niknam\\
Department of Pure Mathematics, Ferdowsi University of Mashhad, Mashhad, P.O. Box 1159-91775, Iran\\

Email: saeedeh.shamsi@gmail.com\\
Email: mjanfada@gmail.com\\
Email: dassamankin@yahoo.co.uk}\\
[2mm]

\end{center}%
\vspace*{0.5cm}
\begin{quotation}
\noindent
{\footnotesize
{\sc Abstract.} In this paper the notion of modular cone metric space is introduced and some properties of such spaces are investigated. Also we define convex modular cone metric which takes values in $C_{\mathbb{R}}(Y)$ where $Y$ is a compact Hausdorff space. Then a fixed point theorem is proved for contractions in these spaces. Furthermore , we make a remark on paper \cite{chwspk} and it will be proved that their fixed point result in modular metric spaces is not true.
}
\end{quotation}

\n \textbf{2010 Mathematics subject classification: } Primary: ; Secondary: .

\markboth
{S. Shamsi Gamchi and M. Janfada, A. Niknam}
 {Modular cone metric spaces}



\section {\sc Introduction and preliminaries}
Ordered normed spaces and cones have applications in applied
mathematics and optimization theory \cite{DN}. Replacing the real
numbers, as the codomain of metrics, by ordered Banach spaces
one may obtain a generalization of metric spaces. Such generalized
spaces called cone metric spaces, were introduced by Rzepecki
\cite{R}. It is proved in \cite{cak} that every cone metric space is metraizable.\\In this paper, which is split into two parts, our aim is to develop the
theory of cone metric spaces called modular cone metric spaces. In the first part, the notion of modular cone metric space is introduced and some properties of such spaces are investigated. Also we define convex modular cone metric which takes values in $C_{\mathbb{R}}(Y)$ where $Y$ is a compact Hausdorff space. Then a fixed point theorem is proved for contractions in these spaces. In the second part of this paper, we make a remark on paper \cite{chwspk} and it will be proved that their fixed point result in modular metric spaces is not true.\\
Let $E$ be a topological vector space (TVS, for short) with its
zero vector $\theta$. A nonempty subset $P$ of $E$ is called a
convex cone if $P+P\subseteq E$ and $\lambda P\subseteq P$ for all
$\lambda\geq 0$. A convex cone $P$ is said to be pointed if
$P\cap(-P)=\{\theta\}$. For a given convex cone $P$ in $E$, a
partial ordering $\preceq$ on $E$ with respect to $P$ is defined
by $x\preceq y$ if and only if $y-x\in P$. We shall write $x\prec y$ if $x\preceq y$ and $x\neq y$, while $x\ll y$ will stand for $y-x\in int P$, where $int P$ denotes the topological interior of $P$.\\
In the sequel, we will need the following useful lemmas:
\begin{lem}\cite{jbr}\label{b}
Let $P$ be a cone in $E$. Then:\\(i) If $\theta \preceq a_n\rightarrow \theta$, then for each $c\in intP$, there exists $N\in\mathbb{N}$ such that for every $n>N$, $a_n\ll c$.\\(ii) For every $c_1,c_2\in intP$, there exists $c\in intP$ such that $c\ll c_1$ and $c\ll c_2$.\\(iii) for every $a\in P$ and $c\in intP$, there exists $n_0\in\mathbb{N}$ such that $a\ll n_0 c$
\end{lem}

The nonlinear scalarization function $\xi_{e}: E\rightarrow
\mathbb{R}$ is
defined as follows:$$\xi_{e}(y)=\inf \{r\in\mathbb{R}:  y\in re-P\}$$for all $y\in E$.\\

\begin{lem}\cite{chy}\label{c}
Let $r\in\mathbb{R}$ and $y\in E$, then
\\(i) $\xi_{e}(y)\leq r
\Leftrightarrow re-y\in P$,\\(ii) $\xi_{e}(y)< r \Leftrightarrow
re-y\in intP$,\\(iii) $\xi_{e}(.)$
 is positively homogeneous and continuous on $E$;\\(iv) If $y_{1}\in y_{2} + P$,
  then $\xi_{e}(y_{2})\leq \xi_{e}(y_{1})$;\\(v) $\xi_{e}(y_{1}+ y_{2})\leq\xi_{e}(y_{1})+\xi_{e}(y_{2})$.
\end{lem}
\begin{defn}\cite{WSD}
Let $X$ be a nonempty set and $d:X\times
X\rightarrow E$ be a mapping that satisfies:\\(CM1) For all $x,y\in X$, $d(x,y)\succeq\theta$ and
$d(x,y)=\theta$ if and only if $x=y$,\\(CM2) $d(x,y)=d(y,x)$ for
all $x,y\in X$,\\(CM3) $d(x,y)\preceq d(x,z)+d(z,y)$  for all
$x,y,z\in X$.\\Then $d$ is called a topological vector space
cone metric (TVS cone metric, for short) on $X$ and $(X,d)$
is said to be a topological vector space cone metric
space.
\end{defn}

\begin{defn}\cite{Vch1}
Let $X$ be a nonempty set. A function $w: (0,\infty)\times X\times X \rightarrow [0,\infty]$ is said to be a modular metric on $X$, if it satisfies the following three axioms:\\
(i) for given $x,y\in X$, $w_{\lambda}(x,y)=0$ for all $\lambda>0$ if and only if $x=y$,\\(ii) $w_\lambda(x,y)=w_\lambda(y,x)$ for all $x,y\in X$ and $\lambda>0$,\\
(iii) $w_{\lambda+\mu}(x,y)\leq w_{\lambda}(x,z)+w_{\mu}(z,y)$ for all $x,y,z\in X$ and $\lambda,\mu\in (0,\infty)$.\\If, instead of (i), the function $w$ satisfies only\\$(i')$ $w_\lambda(x, x) =0$ for all $\lambda>0$ and $x\in X$,\\then $w$ is said to be a pseudo modular on $X$, and if $w$ satisfies $(i')$ and\\
$(i'')$ for given $x,y\in X$, if there exists $\lambda>0$ such that $w_\lambda(x,y)=0$, then $x=y$.\\the function $w$ is called a strict modular on $X$.\\A modular (pseudomodular, strict modular) $w$ on $X$ is said to be convex
if, instead of (iii), for all $\lambda,\mu > 0$ and $x, y, z\in X$ it satisfies the inequality\\$$w_{\lambda+\mu}(x,y)\leq \frac{\lambda}{\lambda+\mu}w_{\lambda}(x,z)+\frac{\mu}{\lambda+\mu}w_{\mu}(z,y)$$
\end{defn}
Note that the conditions (i) to (iii) imply that for all $y,z\in X$ and $\lambda>0$, $w_\lambda(y,z)\geq 0$. If $w_\lambda(x,y)$ does not depend on $x,y\in X$, then by (i) $w\equiv 0$. Now if $w_\lambda(x,y)=w(x,y)$ is independent of $\lambda>0$, then axioms (i)-(iii) mean that $w$ is a metric on X.\\The essential property of a pseudo modular $w$ on $X$ (cf. \cite{Vch1}, Section 2.3)
is that, for any given $x, y\in X$, the function $0 < \lambda\mapsto w_\lambda(x, y)$ is decreasing on $(0,\infty)$.
\begin{defn}\cite{Vch1}
Let $w$ be  pseudo modular on $X$, the two sets
$$X_w=X_w(x_0) = \{x \in X : w_\lambda(x, x_0)\rightarrow 0 ~ as ~ \lambda\rightarrow \infty\}$$
and
$$X_w^*=X_w^*(x_0) = \{x \in X : \exists\lambda=\lambda(x) > 0 ~ \mbox{such that} ~ w_\lambda(x, x_0) < \infty\}$$
are said to be modular spaces (around $x_0$).
\end{defn}
If $w$ is a convex modular on $X$, then according to \cite{Vch1}, Section 3.5 and
Theorem 3.6, the two modular spaces coincide, $X_w = X_w^*$.
\begin{defn}\cite{chwspk}
Let $w$ be a modular metric on $X$, $x\in X_w$ and $\{x_{n}\}$ be a
sequence in $X_w$. Then\\(1) $\{x_{n}\}$ is said to be modular convergent to $x$ if
for every $\lambda>0$, $w_\lambda(x_n,x)\rightarrow 0$ as $n\rightarrow \infty$.\\
(2) $\{x_{n}\}$ is said to be a modular Cauchy sequence if $\lambda>0$, $w_\lambda(x_n,x_m)\rightarrow 0$ as $n,m\rightarrow \infty$\\(3) $X_w$ is called a complete modular metric space if every modular Cauchy sequence is modular convergent.
\end{defn}
\begin{thm}\cite{chist3}\label{d}
Let $w$ be a convex modular metric on $X$ and $X_w^*$ be a complete modular metric space. Suppose that $T: X_w^*\rightarrow X_w^*$ is a mapping which satisfies the following condition:$$w_{k\lambda}(Tx,Ty)\leq w_\lambda(x,y) ~ ~  \mbox{for all} ~ ~ \theta\ll c\ll c_0.$$Then $T$ has a unique fixed point.
\end{thm}

\section {\sc Convergence in modular cone metric spaces}
In the following, unless otherwise specified, we always suppose
that $E$ is a locally convex Hausdorff TVS with its zero vector
$\theta$, $P$ a proper closed and convex pointed cone in $E$ with
$int P\neq\emptyset$, $e\in int P$ and $\preceq$ the partial
ordering with respect to $P$.

Throughout this paper functions $w: intP\times X\times X \rightarrow E$ will be written as $w(c,x,y)=w_c(x,y)$ for all $c\in intP$ and $x,y\in X$.

\begin{defn}\label{a}
Let $X$ be a nonempty set. A function $w: intP\times X\times X \rightarrow E$ is said to be a modular cone metric on $X$, if it satisfies the following
three axioms:\\
(i) For given $x,y\in X$, $w_{c}(x,y)=\theta$ for all $c\in intP$ if and only if $x=y$,\\(ii) $w_c(x,y)=w_c(y,x)$, for all $x,y\in X$ and $c\in intP$,\\
(iii) $w_{c_1+c_2}(x,y)\preceq w_{c_1}(x,z)+w_{c_2}(z,y)$ for all $x,y,z\in X$ and $c_1,c_2\in intP$.\\$(X,E,P,w)$ is called a modular cone metric space.
\end{defn}
Note that the conditions (i) to (iii) imply that for all $y,z\in X$ and $c\in intP$, $\theta\preceq w_c(y,z)$. Indeed, by setting $x=y$ and $c_1=c_2=c$ in (iii), for all $y,z\in X$ one gets $\theta=w_{2c}(y,y)\preceq 2w_c(y,z)$.\\If $w_c(x,y)$ does not depend on $x,y\in X$, then by (i) $w\equiv 0$. Now if $w_c(x,y)=w(x,y)$ is independent of $c\in intP$, then axioms (i)-(iii) mean that $w$ is a cone metric on X.\\
For given $x,y\in X$, the function $\theta\ll c\mapsto w_c(x,y)$ is decreasing on $intP$. In fact, if $\theta\ll c_1\ll c_2$, then $w_{c_2}(x,y)\preceq w_{c_2-c_1}(x,x)+w_{c_1}(x,y)=w_{c_1}(x,y)$.\\
By the following example, we may construct a   of examples of modular cone metric spaces using a cone metric.
\begin{exam}
Let $(X,d)$ be a cone metric space and $\varphi : intP\rightarrow (0,\infty)$ be a decreasing function. Now define $w_c(x,y)=\varphi(c)d(x,y)$. It is easy to see that $w$ is a modular cone metric on $X$. Indeed, by the properties of a cone metric $d$ axioms (i)-(ii) of definition \ref{a} are satisfied. On the other hand, for all $x,y,z\in X$ and $c_1,c_2\in intP$, we have:
\begin{align*}
w_{c_1+c_2}(x,y)=&\varphi(c_1+c_2)d(x,y)\\
&\preceq \varphi(c_1+c_2)d(x,z)+\varphi(c_1+c_2)d(z,y)\\
&\preceq \varphi(c_1)d(x,z)+\varphi(c_2)d(z,y)\\
&=w_{c_1}(x,z)+w_{c_2}(z,y).
\end{align*}
\end{exam}

Now we define some convergence concepts in this space. It will be proved that these definitions are compatible with some topology on $X$ constructed by a modular cone metric.
\begin{defn}
Let $(X,E,P,w)$ be a modular cone metric space, $x\in X$ and $\{x_{n}\}$ be a
sequence in $X$. Then\\(1) $\{x_{n}\}$ is said to be modular cone convergent to $x$ and we denote it by $x_n \overset{w}{\longrightarrow} x$ if
for every $c\gg\theta$ there exists a positive integer $N$ such that for all $n>N$, $w_c(x_{n}, x)\ll c$.\\
(2) $\{x_{n}\}$ is said to be a modular cone Cauchy sequence if for every $c\gg\theta$ there exists a
 positive integer $N$ such that for all $m,n>N$,
$w_c(x_{n}, x_{m})\ll c$.\\(3) $(X, w)$ is called a complete modular cone metric space if every modular cone Cauchy sequence is convergent.
\end{defn}
Let $E$ be a Banach space and $P$ be a cone in $E$. The cone $P$ is called normal if there exists a
constant $K>0$ such that for all $a,b\in P$, $a\preceq b$ implies that $\|a\|\leq K\|b\|$. The least positive number satisfying the above inequality
is called the normal constant of $P$.\\
In the next theorem some equivalent condition for convergence is proved.
\begin{thm}
Let $E$ be a Banach space and $P$ be a normal cone with normal
constant $K$. Suppose that $(X,E,P,w)$ is a modular cone metric space and $\{x_n\}$ is a sequence in $X$. Then,\\
(i) $x_n \overset{w}{\longrightarrow} x$ if and only if for each $c\in intP$, $\|w_c(x_n,x)\|\rightarrow 0$ as $n\rightarrow \infty$.\\
(ii) $\{x_n\}$ is modular cone Cauchy if and only if for each $c\in intP$, $\|w_c(x_n,x_m)\|\rightarrow 0$ as $n,m\rightarrow \infty$.
\end{thm}
\begin{proof}
Let $x_n \overset{w}{\longrightarrow} x$. Fix $c\in intP$. Then for every $0<\varepsilon <1$, there exists a positive integer $N$ such that for all $n>N$, $w_{\varepsilon c}(x_{n}, x)\ll \varepsilon c$. On the other hand, so for each $n$ we have $w_c(x_n,x)\preceq w_{\varepsilon c}(x,y) $, since $\varepsilon c\ll c$. Hence, $w_c(x_n,x)\ll \varepsilon c$, for each $n>N$. Therefore, normality of $P$ implies that $\|w_c(x_n,x)\|\leq K \varepsilon \|c\|$, for each $n>N$. For proving its converse, let for each $c\in intP$, $\|w_c(x_n,x)\|\rightarrow 0$ as $n\rightarrow \infty$. Let $c\in intP$. By assumption we have $w_c(x_n,x)\rightarrow \theta$ as $n\rightarrow \infty$. Applying Lemma \ref{b} there is a positive integer $N$ such that for all $n>N$, $w_c(x_n,x)\ll c$. Thus the proof of (i) is complete.\\For proving (ii) let $\{x_n\}$ be a modular cone Cauchy sequence and fix $c\in intP$. Similar to (i) one can show that $\|w_c(x_n,x_m)\|\rightarrow 0$ as $n,m\rightarrow \infty$, for each $c\in intP$. Conversely, let $w_c(x_n,x_m)\rightarrow \theta$ as $n,m\rightarrow \infty$, for each $c\in intP$. For given $c\in intP$, there exists $\varepsilon>0$ such that $c-y\in intP$ for each $y\in B_{\varepsilon}(\theta)$. For this $\varepsilon>0$, there is a positive integer $N$ such that for all $n,m>N$, $\|w_c(x_n,x_m)\|<\varepsilon$. Therefore, for all $n,m>N$, $c-w_c(x_n,x_m)\in intP$, that means $w_c(x_n,x_m)\ll c$ for all $n,m>N$. Thus the proof is complete.
\end{proof}
Now we are going to construct a topology on $X$, where $(X,E,P,w)$ is a modular cone metric space. For any $x\in X$, $c\in intP$, $B_w(x,c)=\{y\in X: W_c(x,y)\ll c\}$.
\begin{thm}
Let $(X,E,P,w)$ be a modular cone metric space. Then $$\tau_w=\{U\subset X: \forall x\in U ~ \exists c\in intP ~ s.t ~ B_w(x,c)\subset U\}$$forms a Hausdorff topology on $X$.
\end{thm}
\begin{proof}
Trivially $\emptyset , X \in \tau_w$.\\
Also for any $U, V \in \tau_w $ and $x\in U \cap V$, then $x\in U$ and $x \in V$, one may find $c_1, c_2\in intP$ such that $x\in B_w(x, c_1)\subset U$ and
$x\in B_w(x, c_2)\subset V$. By Lemma \ref{b}, there exists $c\in int P$ such that $c\ll c_1$ and $c \ll c_2$. Now suppose that $y\in B_w(x,c)$, so we have $w_{c_1}(x,y)\preceq w_c(x,y)\ll c\ll c_1$, hence, $y\in B_w(x,c_1)$. Similarly we have $y\in B_w(x,c_2)$. Thus $B_w(x, c)\subset B_w(x, c_1)\cap B_w(x, c_2)\subset U\cap V$. Therefore $U\cap V\in \tau_w$.\\
Let $U_\alpha\in \tau_w$ for each $\alpha\in\Delta$ and let $x\in \bigcup_{\alpha\in\Delta} U_{\alpha}$, then $\exists\alpha_{0}\in \Delta$ such that $x\in U_{\alpha_0}$. Hence, $x\in B(x, c)\subset U_{\alpha_0}$, for some $c\in int P$. That is $\bigcup_{\alpha\in\Delta} U_{\alpha}\in\tau$. Thus $\tau_w$ is a topology on $X$.\\
Now we prove that this topology is a Hausdorff topology. Suppose that $x,y\in X$ and $x\neq y$. So by the property $(i)$ of Definition \ref{a} there exists $c\in intP$ such that $w_{c_0}(x,y)\neq \theta$. In contrary, we assume that for each $c\in intP$, $$B_w(x,c)\cap B_w(y,c)\neq\emptyset .$$Hence, for each $n\in\mathbb{N}$, there is a $z_n\in X$ such that $$z_n\in B_w(x,\frac{c_0}{2n})\cap B_w(y,\frac{c_0}{2n}).$$Thus for each $n>1$, we have
\begin{align*}
w_{c_0}(x,y)&\preceq w_{\frac{c_0}{n}}(x,y)\\
&\preceq w_{\frac{c_0}{2n}}(x,z_n)+w_{\frac{c_0}{2n}}(x,z_n)\\
& \ll\frac{c_0}{2n}+\frac{c_0}{2n}\\
&=\frac{c_0}{n}
\end{align*}
Therefore for each $n>1$, $\frac{c_0}{n}-w_{c_0}(x,y)\in P$. But $P$ is closed, so $$\lim_{n\rightarrow\infty}(\frac{c_0}{n}-w_{c_0}(x,y))=-w_{c_0}(x,y)\in P.$$This is a contradiction since $P\cap (-P)=\{\theta\}$. So the proof is complete.
\end{proof}
One can easily see that the collection $\{B_w(x, c): x\in X, c\in intP\}$ forms a basis for $\tau_w$ under which the above definitions of convergent and Cauchy sequences are fully justified.\\
\begin{thm}
Every modular cone metric space $(X,E,P,w)$ is first countable.
\end{thm}
\begin{proof}
Let $x\in X$. Fix $c\in intP$. We show that $\beta_x =\{B_w(x, \frac{c}{n}): n\in\mathbb{N}\}$ is a local base at $x$.
Let $U$ be an open set such that $x\in U$. Therefore, there is a $c_1\in intP$ such that $B_w(x,c_1)\subset U$. on the other hand there exists a positive integer $n$, such that $\frac{c}{n}\ll c_1$. Now suppose that $y\in B_w(x,\frac{c}{n})$, so we have $w_{c_1}(x,y)\preceq w_{\frac{c}{n}}(x,y)\ll \frac{c}{n}\ll c_1$, hence, $y\in B_w(x,c_1)$. Thus $B_w(x, \frac{c}{n})\subset B_w(x, c_1)\subset U$.
\end{proof}
\begin{thm}
Let $(X,E,P,w)$ be a modular cone metric space. A self map $f:X\rightarrow X$ is continuous at $x \in X$ if and only if whenever $x_n \rightarrow x$, we have $f(x_n)\rightarrow f(x)$, as $n\rightarrow \infty$.
\end{thm}
\begin{proof}
Applying Theorem \ref{e} complete the proof.
\end{proof}

The following useful theorem shows that we may construct a family of modular metrics using a modular cone metric and the mapping $\xi_e$.
\begin{thm}\label{e}
Let $(X,E,P,w)$ be a modular cone metric space and $ e\in intP$. Then $W^e: (0,\infty) \times X \times X \rightarrow [0,\infty)$ which is defined by $W_\lambda(x,y)=\xi_e(w_{\lambda e}(x,y))$ is a modular metric on $X$.
\end{thm}
\begin{proof}
Let $x,y\in X$. If $x=y$, then $w_{\lambda e}(x,y)=\theta$ for each $\lambda>0$, so $W_{\lambda}^e(x,y)=\xi_e(w_{\lambda}(x,y))=0$ for each $\lambda>0$. Now, let for each $\lambda>0$, $W_{\lambda}^e(x,y)=0$. So $w_{\lambda e}(x,y)=0\theta$ for each $\lambda>0$. On the other hand for each $c\in intP$ there exists a $\lambda>0$ such that $\lambda e\ll c$, hence $w_c(x,y)\preceq w_{\lambda e}(x,y)=\theta$. Therefore, $w_c(x,y)=\theta $, for each $c\in intP$. That means $x=y$. So we prove that $x=y$ if and only if $W_{\lambda}^e(x,y)=0$, for each $\lambda>0$.\\It is easy to see that $W_{\lambda}^e(x,y)=W_{\lambda}^e(y,x)$, for all $\lambda>0$ and $x,y\in X$.\\Also if $x,y,z\in X$ and $\lambda_1, \lambda_2>0$, then
\begin{align*}
W_{\lambda_1+\lambda_2}^e(x,y)=&\xi_e(w_{\lambda_1 e+\lambda_2 e}(x,y))\\
&\leq \xi_e(w_{\lambda_1 e}(x,z))+\xi_e(w_{\lambda_2 e}(z,y))\\
&= W_{\lambda_1}^e(x,z)+ W_{\lambda_2}^e(z,y)
\end{align*}
\end{proof}

\begin{thm}
Let $(X,E,P,w)$ be a modular cone metric space, $e\in intP$, $x\in X$ and $\{x_n\}$ be a sequence in $X$. Then \\
(i) $x_n \overset{w}{\longrightarrow} x$ if and only if for every $\lambda>0$, $W_{\lambda}^e(x_n,x)\longrightarrow 0$ as $n\longrightarrow \infty$.\\
(ii) $\{x_n\}$ is a  $w$-cauchy sequence if and only if for every $\lambda>0$, $W_{\lambda}^e(x_n,x_m)\longrightarrow 0$ as $n,m\longrightarrow \infty.$
\end{thm}
\begin{proof}
(i) Let $x_n \overset{w}{\longrightarrow} x$. Fix $\lambda>0$. For given $\epsilon >0$ if $\lambda > \epsilon$, then for $c=\epsilon e$, there is a positive integer $N$ such that for each $n>N$, $w_{\epsilon e}(x_n,x)\ll \epsilon e$. On the other hand, so $w_{\lambda e}(x_n,x)\preceq w_{\epsilon e}(x_n,x)\ll \epsilon e$, since $\epsilon e\ll \lambda e$. Hence by the property (iv) of Lemma \ref{c}, for each $n>N$ we have $$W_{\lambda}(x_n,x)\leq W_{\epsilon}(x_n,x)=\xi_e(w_{\epsilon e}(x_n,x))<\xi_e( \epsilon e)=\epsilon.$$Now suppose that $\lambda\leq\epsilon$. For $c=\lambda e$, there is a positive integer $N$ such that for each $n>N$, $w_{\lambda e}(x_n,x)\ll \lambda e$. Hence, for each $n>N$ we have $$W_{\lambda}(x_n,x)=\xi_e(w_{\lambda e}(x_n,x))<\xi_e( \lambda e)=\lambda <\epsilon.$$For its converse, let for each $\lambda >0$, $W_{\lambda}^e(x_n,x)\longrightarrow 0$. Fix $c\in intP$, By Lemma \ref{b}, there exists a positive real number $\lambda_0$ such that $\lambda_0 e\ll c$. So, there is a positive integer $N$ such that for each $n>N$,$$\xi_e(w_{\lambda_0 e}(x_n,x))=W_{\lambda_0}(x_n,x) <\lambda_0 .$$Thus by property (ii) of Lemma \ref{c}  and the fact that the function $\theta\ll c\mapsto w_c(x,y)$ is decreasing on $intP$, we have $$w_{c}(x_n,x)\preceq w_{\lambda_0 e}(x_n,x)\ll \lambda_0 e\ll c,$$for each $n>N$.\\
The proof of (ii) is similar to (i).
\end{proof}

\section{\sc A fixed point theorem for contractions in convex modular cone metric spaces}
In this section we suppose that $E=C_{\mathbb{R}}(Y)$ where $Y$ is a compact Hausdorff space. $E$ with sup-norm is a real Banach space. Put $P=\{f\in C_{\mathbb{R}}(Y): f(x)\geq 0 ~ \mbox{for any} x\in Y
\}$. It is clear that $P$ is a cone in $E$. Partial order on $E$ with respect to $P$ is defined as follows:$$f\preceq g \Leftrightarrow f(x)\leq g(x) ~ \forall x\in Y.$$Note that $intP\neq\emptyset$. Indeed, $1\in intP$.
\begin{lem}
If $f\in intP$, then $f$ is invertible.
\end{lem}

\begin{proof}
Let $f\in intP$. Then there exists $\epsilon >0$ such that $$B_{\epsilon}(f)=\{g\in C_{\mathbb{R}}(Y): \|g-f\|_{\infty}<\epsilon\}\subset P.$$If, in contrary, there is a $t_0 \in Y$ such that $f(t_0)=0$, then by defining $g(t)=f(t)-\frac{\epsilon}{2}$, it is clear that $g\in B_{\epsilon}(f)$ and hence $g\in P$. But $g(t_0)<0$ and this is a contradiction. Thus $0\notin rang(f)$ and so $f$ is invertible.
\end{proof}
Now we give the concept of convex modular cone metric spaces.
\begin{defn}
Let $X$ be a nonempty set. A function $w: intP\times X\times X \rightarrow E$ is said to be a convex modular cone metric on $X$, if it satisfies the following
three axioms:\\
(i) For given $x,y\in X$, $w_{c}(x,y)=\theta$, for all $c\in intP$, if and only if $x=y$,\\(ii) $w_c(x,y)=w_c(y,x)$, for all $x,y\in X$ and $c\in intP$,\\
(iii) $w_{c_1+c_2}(x,y)\preceq \frac{c_1}{c_1+c_2} w_{c_1}(x,z)+\frac{c_2}{c_1+c_2} w_{c_2}(z,y)$, for all $x,y,z\in X$ and $c_1,c_2\in intP$.\\$(X,w)$ is called a convex modular cone metric space.
\end{defn}
\begin{exam}
Let $(X,d)$ be a cone metric space. Now define $w_c(x,y)=\frac{d(x,y)}{c}$. It is easy to see that $w$ is a convex modular cone metric on $X$.
\end{exam}

Note that every convex modular cone metric space is a modular cone metric space.\\
If instead of (i), the function $w$ satisfies in $(i')$ and $(i")$, then it is called strict convex modular cone metric.\\
$(i')$ $w_c(x, x) =\theta$ for all $c\in intP$ and $x\in X$,\\
$(i")$ Given $x,y\in X$, if there exists $
c\in intP$ such that $w_c(x,y)=\theta$, then $x=y$.
\begin{rem}
Let $(X,w)$ be a convex modular cone metric space and $1\in E=C_{\mathbb{R}}(Y)$ be the unit constant function. It is easy to see that $W_{\lambda}^1(x,y)=\xi_1(w_{\lambda 1}(x,y))$ is a convex modular metric on $X$ and $W_{\lambda}^1(x,y)<\infty$ so $X_W^*=X$.
\end{rem}
\begin{thm}
Let $(X,w)$ be a complete convex modular cone metric space. Suppose there exists a constant number $k\in (0,1)$ and $c_0\in intP$ such that a mapping $T: X\rightarrow X$ satisfies the following condition:\\

$$w_{kc}(Tx,Ty)\preceq w_c(x,y) ~ ~  \mbox{for all} ~ ~ \theta\ll c\ll c_0.$$Then $T$ has a unique fixed point.
\end{thm}
\begin{proof}
Take $\lambda_0=\|c_0\|_{\infty}$. For all $0<\lambda <\lambda_0$ by the above inequality for convex modular metric $W_{\lambda}^1(x,y)=\xi_1(w_{\lambda 1}(x,y))$ we have$$W_{k\lambda}^1(x,y)=\xi_1(w_{k\lambda 1}(x,y))\leq \xi_1(w_{\lambda 1}(x,y))=W_{\lambda}^1(x,y).$$Since $X_W^*=X$ and $(X,w)$ is complete, so $(X_{W}^*,W)$ is complete. Hence, the proof is complete by Theorem \ref{d}.
\end{proof}

\section{\sc A note on "Fixed point theorems for contraction mappings
in modular metric spaces"}

Recently, in the paper "Ch. Mongkolkeha, W. Sintunavarat, P. Kumam, Fixed point theorems for contraction mapping in modular metric spaces, Fixed Point Theory and Applications, (2011)", the authors have studied and introduced some fixed-point theorems in
the framework of a modular metric space. We will first state the main result that the authors have proved in that paper and then by constructing an example we will show that the result is not true.
\begin{defn}
(Definition 3.1. of original paper) Let $w$ be a modular metric on $X$ and $X_w$ be a modular metric space
induced by $w$ and $T : X_w \rightarrow X_w$ be an arbitrary mapping. A mapping $T$ is called a contraction
if for each $x, y \in X_w$ and for all $\lambda >0$ there exists $0 < k <1$ such that$$w_\lambda(Tx, Ty)\leq kw_\lambda(x, y).$$
\end{defn}
\begin{thm}
$($Theorem $3.2.$ of original paper$)$ Let $w$ be a modular metric on $X$ and $X_w$ be a modular metric space
induced by $w$. If $X_w$ is a complete modular metric space and $T : X_w \rightarrow X_w$ is a contraction
mapping, then $T$ has a unique fixed point in $X_w$. Moreover, for any $x\in X_w$, iterative
sequence $\{T_nx\}$ converges to the fixed point.
\end{thm}
Now by following example we show that the above theorem is not valid.
\begin{exam}
Let $X = \{(a, 0)\in\mathbb{R}^2: \frac{1}{2}\leq a \leq 1\}\cup \{(0, b) \in\mathbb{R}^2: \frac{1}{2}\leq b \leq 1\}$.
Defined the mapping $w : (0, \infty) \times X \times X\rightarrow [0, \infty]$ by
$$w_\lambda((a_1, 0), (a_2, 0)) =\frac{4|a_1-a_2|}{3\lambda}'$$
$$w_\lambda((0, b_1), (0, b_2)) = \frac{|b_1-b_2|}{\lambda},$$
and $$w\lambda((a, 0), (0, b)) =\frac{4a}{3\lambda}+\frac{b}{\lambda}= w_\lambda((0, b), (a, 0)).$$We note that if we take $\lambda\rightarrow\infty$, then we see that $X = X_w$ and also it is easy to see that if $\{x_n\}$ is a Cauchy sequence in $X_w$ then we just have one of the following assertions:\\$(1)$ There exists a positive integer $N$ such that for each $n>N$, $x_n\in\{(a, 0)\in\mathbb{R}^2: \frac{1}{2}\leq a \leq 1\}$.\\$(2)$ There exists a positive integer $N$ such that for each $n>N$, $x_n\in \{(0, b) \in\mathbb{R}^2: \frac{1}{2}\leq b \leq 1\}$.\\Indeed, It follows from the fact that for $(a,0), (0,b)\in X$, $$w_\lambda((a, 0), (0, b)) =\frac{4a}{3\lambda}+\frac{b}{\lambda}\geq \frac{7}{3\lambda}.$$Thus $X_w$ is a complete modular metric space. Now we define a mapping $T:X_w\rightarrow X_w$ by $$T((a, 0)) = (0, a)$$and$$T((0, b)) =(\frac{b}{2},0)$$Simple computations show that
$$w_\lambda(T((a_1, b_1)), T((a_2, b_2)))\leq \frac{3}{4}w_\lambda((a_1, b_1), (a_2, b_2)),$$for all $(a_1, b_1), (a_2, b_2)\in X_w$. But $T$ does not have any fixed point in $X_w$.
\end{exam}


\providecommand{\bysame}{\leavevmode\hbox
to3em{\hrulefill}\thinspace}


\end{document}